\documentclass[12pt,letterpaper]{amsart}
\usepackage{geometry,mathpazo}
\usepackage{amsmath,amssymb,amsfonts,amsthm,amscd}
\usepackage{mathrsfs,latexsym,url,graphicx,setspace,color,cancel}
\newtheorem{theorem}{Theorem}
\newtheorem*{theorem*}{Theorem}
\newtheorem{proposition}{Proposition}

\theoremstyle{remark}
\newtheorem{remark}{Remark}
\newtheorem{definition}{Definition}

\newcommand{\abs}[1]{\left\lvert#1\right\rvert}
\newcommand{\norm}[1]{\lvert\lvert#1\rvert\rvert}

\newcommand{\disc}{\mathbb{D}}
\newcommand{\C}{\mathbb{C}}
\newcommand{\n}{\mathbb{N}}

\newcommand{\zb}{\overline{z}}
\newcommand{\D}{\Omega}

\newcommand{\dbar}{\overline{\partial}}

\usepackage{hyperref}
\onehalfspace

\title[Hilbert-Schmidt Hankel Operators on Reinhardt Domains]{Hilbert-Schmidt Hankel Operators with anti-holomorphic symbols on Complete Pseudoconvex Reinhardt Domains}

\author{Mehmet \c{C}el\.ik}
\address[Mehmet \c{C}elik]{Texas A\&M University - Commerce, Department of Mathematics,  Binnion Hall, 
Commerce, TX 75429-3011}
\email{mehmet.celik@tamuc.edu}

\author{Yunus E. Zeytuncu}
\address[Yunus E. Zeytuncu]{University of Michigan - Dearborn, Department of Mathematics and Statistics, Dearborn, MI  48128}
\email{zeytuncu@umich.edu}

\subjclass[2000]{Primary 47B35; Secondary 32A36, 47B10}
\keywords{Canonical Solution Operator for $\dbar$-problem, Hankel Operator, Hilbert-Schmidt Operator} 

\begin{document}

\begin{abstract}

On complete pseudoconvex Reinhardt domains in $\C^2$, we show that there is no nonzero Hankel operator with an anti-holomorphic symbol that is Hilbert-Schmidt. We also present examples of  unbounded non-pseudoconvex domains that admit nonzero Hilbert-Schmidt Hankel operators with anti-holomorphic symbols.
\end{abstract}

\maketitle 

\section{Introduction}
\subsection{Setup and Problem}
For a domain $\D$ in $\C^n$, we denote the space of square integrable functions and the space of square integrable holomorphic functions on $\D$ by $L^{2}(\D)$ and $A^{2}(\D)$ (the Bergman space of $\D$), respectively. The Bergman projection operator, $P$, is the orthogonal projection from $L^2(\D)$ onto $A^{2}(\D)$. It is an integral operator with the kernel called the Bergman kernel, which is denoted by $B_{\D}(z,w)$. Moreover, if $\{e_n(z)\}_{n=0}^{\infty}$ is an orthonormal basis for $A^2(\D)$ then the Bergman kernel can be represented as 
$$B_{\D}(z,w)=\sum\limits_{n=0}^{\infty}e_n(z)\overline{e_n(w)}.$$
On complete Reinhardt domains the monomials $\left\{z^{\gamma}\right\}_{\gamma\in \n^{n}}$ (or a subset of them) constitute an orthogonal basis for $A^2(\Omega)$.

For $f\in A^2(\D)$, the Hankel operator with the anti-holomorphic symbol $\overline{f}$ is formally defined on $A^2(\D)$ by
\[ H_{\overline{f}}(g)=(I-P)(\overline{f}g).\]
Note that this (possibly unbounded) operator is densely defined on $A^2(\D)$. 

For a multi-index $\gamma=(\gamma_1,\ldots,\gamma_n)\in \n^{n}$, we set 
\begin{align}
c_{\gamma}^{2}=\int\limits_{\D}\abs{z^{\gamma}}^{2}dV(z).
\end{align}
Then on complete Reinhardt domains the set $\left\{\frac{z^{\gamma}}{c_{\gamma}}\right\}_{\gamma\in \n^{n}}$ gives a complete orthonormal basis for $A^{2}(\D)$. Each $f\in A^{2}(\D)$ can be written in the form $f(z)=\sum\limits_{\gamma\in \n^{n}}f_{\gamma}\frac{z^{\gamma}}{c_{\gamma}}$ where the sum converges in $A^{2}(\D)$, but also uniformly on compact subset of $\D$. For the coefficients $f_{\gamma}$, we have $f_{\gamma}= \langle f(z),\frac{z^{\gamma}}{c_{\gamma}}\rangle_{\D}$.

\begin{definition} 
A linear bounded operator $T$ on a Hilbert space $H$ is called a \textit{Hilbert-Schmidt operator} if there is an orthonormal basis $\{\xi_{j}\}$ for $H$ such that the sum $\sum\limits_{j=1}^{\infty}\norm{T(\xi_{j})}^2$ is finite. 
\end{definition}
The sum does not depend on the choice of orthonormal basis $\{\xi_{j}\}$.  For more on Hilbert-Schmidt operators see \cite[Section X]{Retherford93}.

In this paper, we investigate the following problem. On a given Reinhardt domain in $\C^n$, characterize the symbols for which the corresponding Hankel operators are Hilbert-Schmidt. This question was first studied in $\C$ on the unit disc in \cite{ArazyFisherPeetre88}. The problem was studied on higher dimensional domains in \cite[Theorem at pg. 2]{KeheZhu90} where the author showed that when $n\geq 2$, on an $n$-dimensional complex ball there are no nonzero Hilbert-Schmidt Hankel operators (with anti-holomorphic symbols) on the Bergman space. The result was revisited in \cite{Schneider07} with a more robust approach. On more general domains in higher dimensions, the problem was explored in \cite[Theorem 1.1]{KrantzLiRochberg97} where the authors extended the result \cite[Theorem at pg. 2]{KeheZhu90} to bounded pseudoconvex domains of finite type in $\C^2$ with smooth boundary. Moreover, the authors of the current article studied the same problem on complex ellipsoids \cite{CelikZeytuncu2013}, in $\C^{2}$ with not necessarily smooth boundary.

The same question was investigated on Cartan domains of tube type  in \cite[Section 2]{Arazy1996} and on strongly psuedoconvex domains in \cite{Li93, Peloso94}. Arazy studied the natural generalization of Hankel operators on Cartan domains (a circular, convex, irreducible bounded symmetric domains in $\C^n$) of tube type and rank $r>1$ in $\C^n$ for which $n/r$ is an integer. He showed that there is no non-trivial Hilbert-Schmidt Hankel operators with anti-holomorphic symbols on those type of domains. Li and Peloso, independently, obtained the same result on strongly pseudoconvex domains with smooth boundary. 

\subsection{Results}  Let 
\[\D=\{(z_1,z_2)\in \C^{2}\ |\ z_1\in\disc\ \text{ and }\ |z_2|<e^{-\varphi(z_1)}\}\]
($\varphi(z_1)=\varphi(|z_1|)$) be a complete pseudoconvex Reinhardt domain where monomials $\{z^{\alpha}\}$ (or a subset of monomials) form a complete system for $A^{2}(\D)$. In this paper, we show that on complete pseudoconvex Reinhardt domains in $\C^{2}$, there are no nonzero Hilbert-Schmidt Hankel operators with anti-holomorphic symbols. Moreover, we also present examples of unbounded non-pseudoconvex domains on which there are nonzero Hilbert-Schmidt Hankel operators with anti-holomorphic symbols.

\begin{theorem}\label{Main}
Let $\D$ be as above and $f \in A^2 (\D)$.
If the Hankel operator $H_{\overline{f}}$ is Hilbert-Schmidt on $A^2 (\D)$ then $f$ is constant.
\end{theorem}

\begin{remark}
Theorem \ref{Main} generalizes Zhu's result on the unit ball in $\C^{n}$ \cite{KeheZhu90}, Schnider's result on the unit ball in $\C^{n}$ and its variations \cite{Schneider07}. Theorem \ref{Main} also generalizes the result in \cite[Theorem 1.1]{KrantzLiRochberg97} by dropping the finite type condition on complete pseudoconvex Reinhardt domains.
\end{remark}

\begin{remark}\label{RemarkNew}
The new ingredient in the proof of Theorem \ref{Main} is the explicit use of the pseudoconvexity property of the domain $\D$, see the assumption made at \eqref{assumption 2} and how it is used at \eqref{use of pseudoconvexity}. Additionally, we employ the key estimate \eqref{first step on estimating S_{alpha}} proven in \cite{CelikZeytuncu2013}. 
\end{remark}

\begin{remark}
After completing this note, the authors have learned that by using the estimate \eqref{first step on estimating S_{alpha}}, Le obtained the same result on bounded complete Reinhardt domains without the pseudoconvexity assumption, see \cite{TrieuLe}. Although our statement requires pseudoconvexity, it also works on unbounded domains. The study of complex function theory on unbounded domains (and relation to pseudoconvexity) has been investigated recently in \cite{HarringtonRaich2013,HarringtonRaich2014} and new phenomenas have been observed.  
\end{remark}

Wiegerinck in \cite{Wiegerinck84}, constructed Reinhardt domains (unbounded but with finite volume) in $\C^2$ for which the Bergman spaces are $k$-dimensional. In fact, for these domains the Bergman spaces are spanned by monomials of the form $\{(z_1z_2)^j\}_{j=1}^{k-1}$. Therefore, Hankel operators with non-trivial anti-holomorphic symbols are Hilbert-Schmidt. We revisit these and similar domains in the last section to present examples of domains that admit nonzero Hilbert-Schmidt Hankel operators with anti-holomorphic symbols.

\section{An Identity and an Estimate on Reinhardt Domains}

The set $\left\{\frac{z^{\gamma}}{c_{\gamma}}\right\}_{\gamma\in \n^{n}}$  is an orthonormal basis for $A^{2}(\D)$. In order to prove Theorem \ref{Main}, we will look at the sum 
\begin{align}\label{norm of Hankel - last step}
\sum\limits_{\gamma}\left\|H_{\overline{f}}\left(\frac{z^{\gamma}}{c_{\gamma}}\right)\right\|^2=\sum\limits_{\alpha}\left| f_{\alpha}\right|^{2}\sum\limits_{\gamma}\left(\frac{c_{\alpha+\gamma}^{2}}{c_{\gamma}^{2}}
-\frac{c_{\gamma}^{2}}{c_{\gamma-\alpha}^{2}}\right)
\end{align}
for $f\in A^{2}(\D)$.  For detailed computation of \eqref{norm of Hankel - last step} and of the later estimate \eqref{first step on estimating S_{alpha}} we refer to \cite{CelikZeytuncu2013}.

The term $\sum\limits_{\gamma}\left(\frac{c_{\gamma+\alpha}^{2}}{c_{\gamma}^2}-\frac{c_{\gamma}^{2}}{c_{\gamma-\alpha}^2}\right)$ in the identity \eqref{norm of Hankel - last step} plays an essential role in the rest of the proof, and we label it as,  
\begin{align}\label{S-alpha}
S_{\alpha}:=\sum\limits_{\gamma}\left(\frac{c_{\gamma+\alpha}^{2}}{c_{\gamma}^2}-\frac{c_{\gamma}^{2}}{c_{\gamma-\alpha}^2}\right).
\end{align}
Note that, the Cauchy-Schwarz inequality guarantees that $\frac{c_{\gamma+\alpha}^{2}}{c_{\gamma}^2}-\frac{c_{\gamma}^{2}}{c_{\gamma-\alpha}^2}\geq 0$ for all $\alpha$ and $\gamma$.

The computations above hold on any domains where the monomials (or a subset of monomials) form an orthonormal basis for the Bergman space.

Now, we estimate the term $S_{\alpha}$ on complete pseudoconvex Reinhardt domains. Our goal is to show that $S_{\alpha}$ diverges for all nonzero $\alpha$ on these domains. By \eqref{norm of Hankel - last step}, this will be sufficient to conclude Theorem \ref{Main}. 

In earlier results, $S_{\alpha}$'s were computed explicitly to obtain the divergence. Here we obtain the divergence by using the estimate \eqref{first step on estimating S_{alpha}}.

For any sufficiently large $N$, we have
\begin{align}\label{first step on estimating S_{alpha}}
 S_{\alpha}\geq \sum\limits_{|\gamma|=N} \frac{c_{\gamma+\alpha}^{2}}{c_{\gamma}^2}
\end{align}
for any nonzero $\alpha$, see \cite{CelikZeytuncu2013}.

\section{Computations on Complete Pseudoconvex Reinhardt Domains, proof of Theorem $1$}

Let $\phi(r)\in C^{2}([0,1))$, define the following  complete Reinhardt domain
\[\D=\{(z_1,z_2)\in \C^{2}\ |\ z_1\in\disc\ \text{ and }\ |z_2|<e^{-\phi(z_1)}\}.\]
Note that $\phi(z_1)=\phi(|z_1|)$. 

If $\limsup\limits_{r\rightarrow 1^{-}}\phi(r)$ is finite then $\exists c>0$ such that for any $z_{1}\in\disc$ the fiber in the $z_{2}$ direction contains the disc of radius $c$. Hence, $\D$ contains a polydisc $\disc\times c\disc$. This indicates that there are no nonzero Hilbert-Schmidt Hankel operators with anti-holomorphic symbols on $\D$. 
This also indicates that there are no compact Hankel operators with anti-holomorphic symbols.

Therefore, from  this point we assume  
\[\limsup\limits_{r\rightarrow 1^{-}}\phi(r)=+\infty.\]  
In fact, the later assumption \eqref{assumption 2} made on the domain forces $\phi(r)$ not to oscillate, so we can assume  
\begin{align}\label{assumption 1}
 \lim\limits_{r\rightarrow 1^{-}}\phi(r)=+\infty.
\end{align}
On the other hand, $\D$ is pseudoconvex if and only if $z_1\longrightarrow\phi(|z_1|)$ is a subharmonic function on $\disc$. A simple calculation gives $\Delta\phi(z_1)=\phi^{\prime\prime}(r)+\frac{1}{r}\phi^{\prime}(r)$. We assume $\D$ is pseudoconvex therefore we have 
\begin{align}\label{assumption 2}
 \phi^{\prime\prime}(r)+\frac{1}{r}\phi^{\prime}(r)\geq 0\ \ \text{ on }\ \ (0,1).
\end{align}

Our goal is to show that the sum $\sum\limits_{|\gamma|=N} \frac{c_{\gamma+\alpha}^{2}}{c_{\gamma}^2}$ diverges for any nonzero $\alpha$ on a complete pseudoconvex Reinhardt domain $\D$. We start with computing $c_{\gamma}$'s. 

We have, 
\begin{align*}
c_{\gamma}^2&=\int\limits_{\D}\abs{z^{\gamma}}^{2}dV(z)
=\int\limits_{\disc}\abs{z_1}^{2\gamma_1}\int\limits_{|z_2|<e^{-\phi(|z_1|)}}\abs{z_{2}}^{2\gamma_2}dA(z_{2})dA(z_{1}) \\
&=\int\limits_{\disc}\left\{\abs{z_1}^{2\gamma_1}\frac{2\pi}{2\gamma_2+2} e^{-(2\gamma_{2}+2)\phi(|z_1|)}\right\}dA(z_1)
=\frac{2\pi^{2}}{\gamma_2+1}\int\limits_{0}^{1} r^{2\gamma_1+1}e^{-(2\gamma_{2}+2)\phi(r)}dr.
\end{align*}

For sufficiently large $x$ and $y$, consider the following ratio
\begin{align}\label{the ratio}
R_{x,y}:=
\frac{\int\limits_{0}^{1}r^{x+2\alpha_1} e^{-(y+2\alpha_2)\phi(r)}dr}
{\int\limits_{0}^{1}r^{x}e^{-y\phi(r)}dr},
\end{align}
and define 
\[\Phi_{x,y}(r):=
\frac{r^{x}e^{-y\phi(r)}}
{\int\limits_{0}^{1}r^{x}e^{-y\phi(r)}dr}.\]
Note that $\Phi_{x,y}(0)=0$, $\Phi_{x,y}(1)=0$, and $\int\limits_{0}^{1}\Phi_{x,y}(r)dr=1$.

Also, define 
\begin{align}\label{g-alpha}
g_{\alpha}(r)=r^{2\alpha_1}e^{-2\alpha_2\phi(r)}.
\end{align}
Note that $g_{\alpha}(r)$ does not vanish inside the interval $(0,1)$, but may vanish at $r=0$ and $r=1$ depending on $\alpha$. Now, we can rewrite the ratio $R_{x,y}$ as 
\begin{align}
 R_{x,y}=
\int\limits_{0}^{1}\Phi_{x,y}(r)r^{2\alpha_1}e^{-2\alpha_2\phi(r)}dr
=\int\limits_{0}^{1}\Phi_{x,y}(r)g_{\alpha}(r)dr.
\end{align}

Our goal is to find a sub-interval $(a,b)\subset\subset(0,1)$ such that for sufficiently large $x$ and $y$
\[\int\limits_{a}^{b}\Phi_{x,y}(r)dr\geq \frac{1}{2}.\]

For this purpose, we analyze $\Phi_{x,y}(r)$ further on $(0,1)$ and locate the local maximum of $\Phi_{x,y}(r)$. We have
\[\frac{d}{dr}\Phi_{x,y}(r)=\left(x-y\phi^{\prime}(r)r\right)\left(r^{x-1}e^{-y\phi(r)}\right)\frac{1}
{\left(\int\limits_{0}^{1}r^{x}e^{-y\phi(r)}dr\right)}.\]
Therefore, 
\[\frac{d}{dr}\Phi_{x,y}(r)=0\ \text{ on }\ (0,1)\ \text{ when }\ x-y\phi^{\prime}(r)r=0.\]
We label $f_{x,y}(r):=x-y\phi^{\prime}(r)r$. Note that $f_{x,y}(r)$ controls the sign of $\frac{d}{dr}\Phi_{x,y}(r)$, since the rest of the terms in $\frac{d}{dr}\Phi_{x,y}(r)$ are positive. Furthermore, 
\[f_{x,y}(0)=x\ >0\] and 
\begin{align}\label{use of pseudoconvexity}
\frac{d}{dr}f_{x,y}(r)=-y\left(\phi^{\prime}(r)+r\phi^{\prime\prime}(r)\right)\ <0\ \text{( by the assumption \eqref{assumption 2}}).
\end{align}
Hence, $f_{x,y}(r)$ decreases on $(0,1)$ and can vanish at a point. We will show that by choosing $x,y$ appropriately we can guarantee $f_{x,y}(r)$ vanishes on $(0,1)$.\\
All we need is a point $s\in(0,1)$ such that 
\[s\phi^{\prime}(s)>0.\]
However, this is possible by the assumption \eqref{assumption 1}. If there was no such a point $s\in(0,1)$, then $\phi(r)$ wouldn't grow up to infinity. Moreover, if $\exists s\in(0,1)$ such that $s\phi^{\prime}(s)>0$ then since $r\phi^{\prime}(r)>0$ is an increasing function 
\[r\phi^{\prime}(r)>0\ \text{ for all }\ r\in[s,1).\]
 Therefore, there exists a relatively compact subinterval $(a,b)$ of $(0,1)$ such that
\[a\phi^{\prime}(a)>0\]
and hence $r\phi^{\prime}(r)>0$ on $(a,b)$. Moreover, by choosing $x$ and $y$ appropriately we can make 
\begin{align*}
 f_{x,y}\left(a\right)>0 ~\text{ and }~ f_{x,y}\left(b\right)<0.
\end{align*}
That is,
\[x-ya\phi^{\prime}(a)>0  ~\text{  and }~  x-yb\phi^{\prime}(b)<0.\]
Equivalently,
\[a\phi^{\prime}(a)<\frac{x}{y} ~\text{ and }~ \frac{x}{y}<b\phi^{\prime}(b).\]
Therefore, as long as we keep 
\begin{align}\label{the interval}
a\phi^{\prime}(a)<\frac{x}{y}<b\phi^{\prime}(b)
\end{align}
there exist a solution to $x-yr\phi^{\prime}(r)=0$ on the interval $(a,b)\subset\subset(0,1)$, and so we guarantee that the function $\Phi_{x,y}(r)$ assumes its maximum somewhere inside $(a,b)$. \\
Let us take the point $\rho_{xy}\in(a,b)$ where $\Phi_{x,y}(r)$ takes its maximum value . We have 
\[\int\limits_{0}^{\frac{a}{2}}\Phi_{x,y}(r)dr\leq \int\limits_{\frac{a}{2}}^{\rho_{xy}}\Phi_{x,y}(r)dr\ \ \text{ and }\ \ \int\limits_{\frac{1+b}{2}}^{1}\Phi_{x,y}(r)dr\leq \int\limits_{\rho_{xy}}^{\frac{1+b}{2}}\Phi_{x,y}(r)dr\]
Hence, we deduce that
\begin{align}\label{key inequality} 
\int\limits_{\frac{a}{2}}^{\frac{1+b}{2}}\Phi_{x,y}(r)dr\geq \int\limits_{0}^{1}\Phi_{x,y}(r)dr\geq \frac{1}{2}
\end{align}
as long as $a\phi^{\prime}(a)<\frac{x}{y}<b\phi^{\prime}(b)$.
The inequality at \eqref{key inequality} is the crucial step for the rest of the proof. It guarantees that the integral of $\Phi_{x,y}(r)$ is located somewhere in the middle, i.e. does not lean towards any of the end points.

For a multi-index $\gamma=(\gamma_1,\gamma_2)$, let us write $\Phi_{\gamma}(r)=\Phi_{\gamma_1,\gamma_2}(r)$. Then
\begin{align}
 \frac{c_{\gamma+\alpha}^2}{c_{\gamma}^2}&=\frac{\gamma_2+1}{\gamma_2+\alpha_2+1}
\cdot\frac{\int\limits_{0}^{1}r^{2\gamma_1+2\alpha_1+1} e^{-(2\gamma_2+2+2\alpha_2)\phi(r)}dr}
{\int\limits_{0}^{1}r^{2\gamma_1+1}e^{-(2\gamma_2+2)\phi(r)}dr}\\
\nonumber&=\frac{\gamma_2+1}{\gamma_2+\alpha_2+1}\int\limits_{0}^{1}\Phi_{2\gamma_1+1,2\gamma_2+2}(r)g_{\alpha}(r)dr
\end{align}
Then, 
\begin{align}\label{First step on the estimate}
 S_{\alpha}&\geq \sum\limits_{|\gamma|=N} \frac{c_{\gamma+\alpha}^{2}}{c_{\gamma}^2}
=\sum\limits_{k=0}^{N} \frac{c_{\alpha+(k,N-k)}^{2}}{c_{(k,N-k)}^2}\\
\nonumber &=\sum\limits_{k=0}^{N} \frac{c_{(k+\alpha_1,N-k+\alpha_2)}^{2}}{c_{(k,N-k)}^2}=\sum\limits_{k=0}^{N}\frac{N-k+1}{N-k+\alpha_2+1}\int\limits_{0}^{1}\Phi_{2k+1,2(N-k)+2}(r)g_{\alpha}(r)dr.\\
\end{align}

We want to keep \[\frac{2k+1}{2N-2k+2}\in\left(a\phi^{\prime}(a),b\phi^{\prime}(b)\right),\] see \eqref{the interval}. It is equivalent to asking $k$ to be in the interval
\begin{align*}
 \frac{2a\phi^{\prime}(a)}{2a\phi^{\prime}(a)+2}N+\frac{2a\phi^{\prime}(a)-1}{2a\phi^{\prime}(a)+2}<k<\frac{2b\phi^{\prime}(b)}{2b\phi^{\prime}(b)+2}N+\frac{2b\phi^{\prime}(b)-1}{2b\phi^{\prime}(b)+2}.
\end{align*}
We further restrict $k$ to the interval 
\begin{align*}
I_{N}:= \left(\frac{2a\phi^{\prime}(a)}{2a\phi^{\prime}(a)+2}N+\frac{2a\phi^{\prime}(a)-1}{2a\phi^{\prime}(a)+2}\ ,\ \frac{2b\phi^{\prime}(b)}{2b\phi^{\prime}(b)+2}N+\frac{2b\phi^{\prime}(b)-1}{2b\phi^{\prime}(b)+2}\right)\cap (0,N).
\end{align*}
Therefore, the estimate \eqref{First step on the estimate} can be rewritten as
\begin{align}
 S_{\alpha}\geq \sum\limits_{k\in I_{N}}\frac{N-k+1}{N-k+\alpha_2+1}\int\limits_{0}^{1}\Phi_{2k+1,2(N-k)+2}(r)g_{\alpha}(r)dr.
\end{align}
When $k\in I_{N}$ we have
\begin{align*}
\frac{N-k+1}{N-k+\alpha_2+1}\int\limits_{0}^{1}\Phi_{2k+1,2(N-k)+2}(r)g_{\alpha}(r)dr&\geq \frac{1}{1+\alpha_2}\int\limits_{\frac{a}{2}}^{\frac{1+b}{2}}\Phi_{2k+1,2(N-k)+2}(r)g_{\alpha}(r)dr\\
&\geq\frac{1}{1+\alpha_2}\left(\min\limits_{\frac{a}{2}\leq r\leq\frac{1+b}{2}}\{g_{\alpha}(r)\}\right)\int\limits_{\frac{a}{2}}^{\frac{1+b}{2}}\Phi_{2k+1,2(N-k)+2}(r)dr\\ 
\text{by \eqref{key inequality} }\ \ \ &\geq\frac{1}{1+\alpha_2}\left(\min\limits_{\frac{a}{2}\leq r\leq\frac{1+b}{2}}\{g_{\alpha}(r)\}\right)\frac{1}{2}.
\end{align*}
Let $\lambda_{\alpha}:=\frac{1}{2(1+\alpha_2)}\left(\min\limits_{\frac{a}{2}\leq r\leq\frac{1+b}{2}}\{g_{\alpha}(r)\}\right)$. Note that $\lambda_{\alpha}>0$ since $g_{\alpha}(r)$ is strictly positive on $\left(\frac{a}{2},\frac{1+b}{2}\right)$, see \eqref{g-alpha}. This gives us
\begin{align*}
 S_{\alpha}&\geq \sum\limits_{k\in I_N} \frac{c_{\gamma+\alpha}^{2}}{c_{\gamma}^2}\geq\sum\limits_{k\in I_{N}}\frac{N-k+1}{N-k+\alpha_2+1}\int\limits_{0}^{1}\Phi_{2k+1,2(N-k)+2}(r)g_{\alpha}(r)dr
\geq \sum\limits_{k\in I_N} \lambda_{\alpha}=|I_{N}|\lambda_{\alpha}.
\end{align*}
Note that the number of integers in $I_{N}$ is comparable to $N$. Therefore, $S_{\alpha}\gtrsim N$ and this suffices to conclude $S_{\alpha}$ diverges for nonzero $\alpha$.

\section{Examples of  Unbounded Non-Pseudoconvexs Domain with Nonzero Hilbert-Schmidt Hankel Operators}

In this section, we present two examples of domains that admit nonzero Hilbert-Schmidt Hankel operators with anti-holomorphic symbols. In the first example, the Bergman space is finite dimensional and the claim holds for trivial reasons. In the second example, the Bergman space is infinite dimensional; however, some of the terms $S_{\alpha}$'s are bounded. 

We start with defining the following domains from \cite{Wiegerinck84}.
\begin{align*}
X_{1}&=\left\{(z_1,z_2)\in\C^2\ : |z_1|>e, |z_2|<\frac{1}{|z_1|\log|z_1|}\right\}\\
X_{2}&=\left\{(z_1,z_2)\in\C^2\ : |z_2|>e, |z_1|<\frac{1}{|z_2|\log|z_2|}\right\}\\
X_{3}&=\left\{(z_1,z_2)\in\C^2\ : |z_1|\leq e, |z_2|\leq e\right\}\\
\D_0&=X_1\cup X_2\cup X_3\\
B_m&=\left\{(z_1,z_2)\in\C^2\ :|z_1|, |z_2|>1,\Bigl\lvert |z_1|- |z_2|\Bigr\rvert<\frac{1}{(|z_1|+|z_2|)^m}\right\}\\
\D_k&=\D_0\cup B_{4k}
\end{align*}
Note that $\D_0$ and $\D_k$ are unbounded non-pseudoconvex complete Reinhardt domains with finite volume. The following proposition is also from \cite{Wiegerinck84}.
\begin{proposition}\label{Proposition 1} Let $k$ be a positive integer.
\begin{itemize}
\item[(i.)] The Bergman space, $A^2(\D_k)$, is spanned by the monomials $\left\{(z_1z_2)^j\right\}_{j=0}^k$.
\item[(ii.)] The Bergman space, $A^2(\D_0)$, is spanned by the monomials $\left\{(z_1z_2)^j\right\}_{j=0}^{\infty}$.
\end{itemize}
\end{proposition}

Next, we look at the Hankel operators on the Bergman spaces of $\D_0$ and $\D_k$.

\subsection{Example 1} We start with $\D_k$. Since $A^2(\D_k)$ is finite dimensional, for any multi-index of the form $(j,j)$ for $j=1,\cdots,k$; the term $S_{(j,j)}$ is a finite sum and consequently finite when restricted on the subspace of $A^2(\D_k)$ where the multiplication operator with the symbol $\overline{f}$ is bounded. Hence, for any $f\in A^2(\D_k)$, the Hankel operator with the symbol $\overline{f}$ is Hilbert-Schmidt on the subspace of $A^2(\D_k)$ where the operator is bounded.


\subsection{Example 2} Next, we look at $\D_0$ and we observe that the terms $S_{\alpha}$ take a simpler form. Namely, for a multi-index $(j,j)$, 
$$S_{(j,j)}=\sum_{k=0}^{\infty}\left(\frac{c_{(k+j,k+j)}^2}{c_{(k,k)}^2}-\frac{c_{(k,k)}^2}{c_{(k-j,k-j)}^2}\right),$$
where
$$c_{(k,k)}^2=\int_{\D_0}|z_1z_2|^{2k}dV(z_1,z_2).$$

We will particularly compute $S_{(1,1)}$. A simple integration indicates, 
$$c_{(k,k)}^2=4\pi^2\left(\frac{2}{2k+1}+\frac{e^{4k+4}}{(2k+2)^2}\right)$$
and with simple algebra we obtain,
\begin{align*}
\frac{c_{(k+1,k+1)}^2}{c_{(k,k)}^2}-\frac{c_{(k,k)}^2}{c_{(k-1,k-1)}^2}&=\frac{e^{8k+8}\frac{(2k+2)^4-(2k+4)^2(2k)^2}{(2k+4)^2(2k)^2(2k+2)^4}+e^{4k}\frac{p_1(k)}{p_2(k)}+\frac{p_3(k)}{p_4(k)}}{e^{8k+8}\frac{1}{(2k)^2(2k+2)^2}+e^{4k}\frac{p_5(k)}{p_6(k)}+\frac{p_7(k)}{p_8(k)}}
\end{align*}
where $p_1(k),\cdots,p_8(k)$ are polynomials in $k$. For large values of $k$, the first terms at the numerator and the denominator dominate and we obtain,
$$\frac{c_{(k+1,k+1)}^2}{c_{(k,k)}^2}-\frac{c_{(k,k)}^2}{c_{(k-1,k-1)}^2}\approx \frac{\frac{(2k+2)^4-(2k+4)^2(2k)^2}{(2k+4)^2(2k)^2(2k+2)^4}}{\frac{1}{(2k)^2(2k+2)^2}}\approx \frac{1}{k^2}.
$$
Therefore, $S_{(1,1)}$ is finite and the Hankel operator $H_{\overline{z_1z_2}}$ is Hilbert-Schmidt on $A^2(\D_0)$. 

\section{Remarks}
\subsection{Canonical solution operator for $\dbar$-problem:}
The canonical solution operator for $\dbar$-problem restricted to $(0,1)$-forms with holomorphic coefficients is not a Hilbert-Schmidt operator on complete pseudoconvex Reinhardt domains because the canonical solution operator for $\dbar$-problem restricted to $(0,1)$-forms with holomorphic coefficients is a sum of Hankel operators with $\{\zb_{j}\}_{j=1}^{n}$ as symbols (by Theorem \ref{Main} such Hankel operators are not Hilbert-Schmidt), $$\dbar^*N_1(g)=\dbar^*N_1 \left(\sum\limits_{j=1}^{n}g_{j}d\zb_{j}\right)=\sum\limits_{j=1}^{n}H_{\zb_j}(g_{j})$$ 
for any $(0,1)$-form $g$ with holomorphic coefficients. 

\section{Acknowledgement}
We would like to thank Trieu Le for valuable comments on an earlier version of this manuscript.

\bibliographystyle{amsplain}
\bibliography{CelikZeytuncuBib}
\end{document}